\documentclass[11pt]{amsart}
\usepackage{amsfonts}
\usepackage{amsmath}
\usepackage{amssymb}
\usepackage{wasysym}
\usepackage{txfonts}
\usepackage{color}
\usepackage{graphicx}
\usepackage{xspace}
\usepackage{axodraw}
 \font \eightrm=cmr8

 \newcommand{\nc}{\newcommand}

 \nc{\butcher}{{\scriptstyle\ \circleright\ }}

\nc{\surj}{\to\hskip -3mm \to}

\def\racine{{\scalebox{0.3}{
\begin{picture}(12,12)(38,-38)
\SetWidth{0.5} \SetColor{Black} \Vertex(45,-30){6}
\end{picture}
}}}
  
 \def\arbrea{\,{\scalebox{0.15}{ 
  \begin{picture}(8,55) (370,-248)
    \SetWidth{2}
    \SetColor{Black}
    \Line(374,-244)(374,-200)
    \Vertex(374,-197){9}
    \Vertex(375,-245){12}
  \end{picture}
}}\,}

 \def\arbreba{\,{\scalebox{0.15}{ 
\begin{picture}(8,106) (370,-197)
    \SetWidth{2}
    \SetColor{Black}
    \Line(374,-193)(374,-149)
    \Vertex(374,-146){9}
    \Vertex(375,-194){12}
    \Line(374,-142)(374,-98)
    \Vertex(374,-95){9}
  \end{picture}
}}\,}

 \def\arbrebb{\,{\scalebox{0.15}{ 
  \begin{picture}(48,48) (349,-255)
    \SetWidth{2}
    \SetColor{Black}
    \Vertex(375,-252){12}
    \Line(376,-250)(395,-215)
    \Line(373,-251)(354,-214)
    \Vertex(353,-211){9}
    \Vertex(395,-213){9}
  \end{picture}
}}}

\def\arbreca{\,{\scalebox{0.15}{
\begin{picture}(8,156) (370,-147)
    \SetWidth{2}
    \SetColor{Black}
    \Line(374,-143)(374,-99)
    \Vertex(374,-96){9}
    \Vertex(375,-144){12}
    \Line(374,-92)(374,-48)
    \Vertex(374,-45){9}
    \Line(374,-42)(374,2)
    \Vertex(374,5){9}
  \end{picture}
}}\,}

\def\arbrecaB{\,{\scalebox{0.25}{
\begin{picture}(8,156) (370,-147)
    \SetWidth{2}
    \SetColor{Black}
    \Line(374,-143)(374,-99)
    \Vertex(374,-96){9}
    \Vertex(375,-144){12}
    \Line(374,-92)(374,-48)
    \Vertex(374,-45){9}
    \Line(374,-42)(374,2)
    \Vertex(374,5){9}
  \end{picture}
}}\,}

\def\arbrecb{\,{\scalebox{0.15}{
\begin{picture}(48,94) (349,-255)
\SetWidth{2}
\SetColor{Black}
\Line(376,-204)(395,-169)
\Line(373,-205)(354,-168)
\Vertex(353,-165){9}
\Vertex(395,-167){9}
\Vertex(374,-205){9}
\Line(374,-246)(374,-209)
\Vertex(374,-252){12}
\end{picture}}}\,}

\def\arbrecbB{\,{\scalebox{0.25}{
\begin{picture}(48,94) (349,-255)
\SetWidth{2}
\SetColor{Black}
\Line(376,-204)(395,-169)
\Line(373,-205)(354,-168)
\Vertex(353,-165){9}
\Vertex(395,-167){9}
\Vertex(374,-205){9}
\Line(374,-246)(374,-209)
\Vertex(374,-252){12}
\end{picture}}}\,}

\def\arbrecc{\,{\scalebox{0.15}{
 \begin{picture}(48,98) (349,-205)
    \SetWidth{2}
    \SetColor{Black}
    \Vertex(375,-202){12}
    \Line(376,-200)(395,-165)
    \Line(373,-201)(354,-164)
    \Vertex(353,-161){9}
    \Vertex(395,-163){9}
    \Line(353,-160)(353,-113)
    \Vertex(353,-111){9}
  \end{picture}
}}\,}

\def\arbreccB{\,{\scalebox{0.25}{
 \begin{picture}(48,98) (349,-205)
    \SetWidth{2}
    \SetColor{Black}
    \Vertex(375,-202){12}
    \Line(376,-200)(395,-165)
    \Line(373,-201)(354,-164)
    \Vertex(353,-161){9}
    \Vertex(395,-163){9}
    \Line(353,-160)(353,-113)
    \Vertex(353,-111){9}
  \end{picture}
}}\,}

\def\arbrecd{\,{\scalebox{0.15}{
\begin{picture}(48,52) (349,-251)
    \SetWidth{2}
    \SetColor{Black}
    \Vertex(375,-248){12}
    \Line(376,-246)(395,-211)
    \Line(373,-247)(354,-210)
    \Vertex(353,-207){9}
    \Vertex(395,-209){9}
    \Line(375,-247)(375,-206)
    \Vertex(376,-203){9}
  \end{picture}
 }}\,}
 
 \def\arbrecdB{\,{\scalebox{0.25}{
\begin{picture}(48,52) (349,-251)
    \SetWidth{2}
    \SetColor{Black}
    \Vertex(375,-248){12}
    \Line(376,-246)(395,-211)
    \Line(373,-247)(354,-210)
    \Vertex(353,-207){9}
    \Vertex(395,-209){9}
    \Line(375,-247)(375,-206)
    \Vertex(376,-203){9}
  \end{picture}
 }}\,}

\def\arbrece{\,{\scalebox{0.2}{
  \begin{picture}(45,78) (7,-37)
    \SetWidth{1.0}
    \SetColor{Black}
    \Line(30,-27)(15,3)
    \SetWidth{0.5}
    \Vertex(30,-27){9.9}
    \SetWidth{1.0}
    \Line(30,-27)(45,3)
    \SetWidth{0.5}
    \Vertex(15,3){7.07}
    \Vertex(45,3){7.07}
    \SetWidth{1.0}
    \Line(45,3)(45,33)
    \SetWidth{0.5}
    \Vertex(45,33){7.07}
  \end{picture}
}}\,}

\def\arbreceB{\,{\scalebox{0.33}{
  \begin{picture}(45,78) (7,-37)
    \SetWidth{1.0}
    \SetColor{Black}
    \Line(30,-27)(15,3)
    \SetWidth{0.5}
    \Vertex(30,-27){9.9}
    \SetWidth{1.0}
    \Line(30,-27)(45,3)
    \SetWidth{0.5}
    \Vertex(15,3){7.07}
    \Vertex(45,3){7.07}
    \SetWidth{1.0}
    \Line(45,3)(45,33)
    \SetWidth{0.5}
    \Vertex(45,33){7.07}
  \end{picture}
}}\,}

\def\arbreda{\,{\scalebox{0.15}{
\begin{picture}(8,204) (370,-99)
    \SetWidth{2}
    \SetColor{Black}
    \Line(374,-95)(374,-51)
    \Vertex(374,-48){9}
    \Vertex(375,-96){12}
    \Line(374,-44)(374,0)
    \Vertex(374,3){9}
    \Line(374,6)(374,50)
    \Vertex(374,53){9}
    \Line(374,53)(374,98)
    \Vertex(374,101){9}
  \end{picture}
}}\,}

\def\arbredb{\,{\scalebox{0.15}{
\begin{picture}(48,135) (349,-255)
    \SetWidth{2}
    \SetColor{Black}
    \Line(376,-163)(395,-128)
    \Line(373,-164)(354,-127)
    \Vertex(353,-124){9}
    \Vertex(395,-126){9}
    \Vertex(374,-164){9}
    \Line(374,-205)(374,-168)
    \Vertex(374,-207){9}
    \Line(374,-248)(374,-211)
    \Vertex(374,-252){12}
  \end{picture}
}}\,}

\def\arbredc{\,{\scalebox{0.15}{
 \begin{picture}(48,150) (349,-205)
    \SetWidth{2}
    \SetColor{Black}
    \Line(376,-148)(395,-113)
    \Line(373,-149)(354,-112)
    \Vertex(353,-109){9}
    \Vertex(395,-111){9}
    \Line(353,-108)(353,-61)
    \Vertex(353,-59){9}
    \Line(374,-200)(374,-153)
    \Vertex(374,-149){9}
    \Vertex(374,-202){12}
  \end{picture}
}}\,}

\def\arbredd{\,{\scalebox{0.15}{
 \begin{picture}(48,99) (349,-251)
    \SetWidth{2}
    \SetColor{Black}
    \Line(376,-199)(395,-164)
    \Line(373,-200)(354,-163)
    \Vertex(353,-160){9}
    \Vertex(395,-162){9}
    \Vertex(376,-156){9}
    \Vertex(376,-248){12}
    \Line(375,-245)(375,-204)
    \Line(375,-200)(375,-159)
    \Vertex(375,-201){9}
  \end{picture}
}}\,}

\def\arbrede{\,{\scalebox{0.15}{
 \begin{picture}(48,153) (349,-150)
    \SetWidth{2}
    \SetColor{Black}
    \Vertex(375,-147){12}
    \Line(376,-145)(395,-110)
    \Line(373,-146)(354,-109)
    \Vertex(353,-106){9}
    \Vertex(395,-108){9}
    \Line(353,-105)(353,-58)
    \Vertex(353,-56){9}
    \Line(353,-52)(353,-5)
    \Vertex(353,-1){9}
  \end{picture}
}}\,}

\def\arbredf{\,{\scalebox{0.15}{
\begin{picture}(48,98) (349,-205)
    \SetWidth{2}
    \SetColor{Black}
    \Vertex(375,-202){12}
    \Line(376,-200)(395,-165)
    \Line(373,-201)(354,-164)
    \Vertex(353,-161){9}
    \Vertex(395,-163){9}
    \Line(353,-160)(353,-113)
    \Vertex(353,-111){9}
    \Line(395,-159)(395,-112)
    \Vertex(395,-111){9}
  \end{picture}
}}\,}

\def\arbredz{\,{\scalebox{0.15}{
  \begin{picture}(68,88) (329,-215)
    \SetWidth{2}
    \SetColor{Black}
    \Vertex(375,-212){12}
    \Line(376,-210)(395,-175)
    \Line(373,-211)(354,-174)
    \Vertex(353,-171){9}
    \Vertex(395,-173){9}
    \Line(351,-168)(332,-131)
    \Line(355,-168)(374,-133)
    \Vertex(333,-131){9}
    \Vertex(374,-131){9}
  \end{picture}
}}\,}

\def\arbredg{\,{\scalebox{0.15}{
\begin{picture}(48,98) (349,-205)
    \SetWidth{2}
    \SetColor{Black}
    \Vertex(375,-202){12}
    \Line(376,-200)(395,-165)
    \Line(373,-201)(354,-164)
    \Vertex(353,-161){9}
    \Vertex(395,-163){9}
    \Line(375,-201)(375,-160)
    \Vertex(376,-157){9}
    \Vertex(376,-111){9}
    \Line(375,-155)(375,-114)
  \end{picture}
}}\,}

\def\arbredh{\,{\scalebox{0.15}{
 \begin{picture}(90,46) (330,-257)
    \SetWidth{2}
    \SetColor{Black}
    \Vertex(375,-254){12}
    \Line(376,-252)(395,-217)
    \Vertex(395,-215){9}
    \Line(374,-254)(335,-226)
    \Vertex(334,-224){9}
    \Line(375,-252)(356,-215)
    \Vertex(355,-215){9}
    \Line(374,-255)(417,-227)
    \Vertex(418,-225){9}
  \end{picture}
}}\,}

\def\hab{\,{\scalebox{0.25}{
   \begin{picture}(88,101) (314,-236)
    \SetWidth{3}
    \SetColor{Black}
    \GOval(360,-184)(40,40)(0){0.882}
    \Vertex(394,-169){8.49}
    \Vertex(324,-170){8.49}
    \Vertex(361,-226){9.9}
    \end{picture}}}\,}  

\def\hac{\,{\scalebox{0.35}{
   \begin{picture}(88,101) (314,-236)
    \SetWidth{3}
    \SetColor{Black}
    \GOval(360,-184)(40,40)(0){0.882}
    \Vertex(394,-169){8.49}
    \Vertex(324,-170){8.49}
    \Vertex(361,-226){9.9}
    \Vertex(358,-145){8.49}
  \end{picture}}}\,}
  
  \def\hacab{\,{\scalebox{0.25}{
  \begin{picture}(88,171) (314,-245)
    \SetWidth{3}
    \SetColor{Black}
    \GOval(360,-114)(40,39)(0){0.882}
    \Vertex(362,-235){9.9}
    \Line(362,-234)(362,-146)
    \Vertex(361,-154){8.49}
    \Vertex(394,-99){8.49}
    \Vertex(324,-100){8.49}
  \end{picture}}}\,}
  
   \def\hacbe{\,{\scalebox{0.25}{
  \begin{picture}(134,147) (314,-190)
    \SetWidth{3}
    \SetColor{Black}
    \GOval(360,-138)(40,39)(0){0.882}
    \Vertex(394,-123){8.49}
    \Vertex(324,-124){8.49}
    \Line(396,-119)(437,-56)
    \Vertex(440,-53){8.49}
    \Vertex(361,-180){9.9}
  \end{picture}}}\,}
  
  \def\hacca{\,{\scalebox{0.25}{
  \begin{picture}(137,137) (265,-200)
    \SetWidth{3}
    \SetColor{Black}
    \GOval(360,-148)(40,39)(0){0.882}
    \Vertex(394,-133){8.49}
    \Vertex(324,-134){8.49}
    \Vertex(361,-190){9.9}
    \Line(325,-136)(277,-78)
    \Vertex(275,-73){8.49}
  \end{picture}}}\,}
 
  \def\hacdc{\,{\scalebox{0.25}{ 
\begin{picture}(143,126) (314,-245)
    \SetWidth{3}
    \SetColor{Black}
    \Line(361,-234)(361,-194)
    \GOval(360,-159)(40,39)(0){0.882}
    \Vertex(394,-144){8.49}
    \Vertex(361,-235){9.9}
    \Line(363,-234)(451,-178)
    \Vertex(449,-176){8.49}
    \Vertex(324,-145){8.49}
  \end{picture}}}\,}
  
  \def\haced{\,{\scalebox{0.25}{ 
  \begin{picture}(142,126) (260,-245)
    \SetWidth{3}
    \SetColor{Black}
    \Line(361,-234)(361,-194)
    \GOval(360,-159)(40,39)(0){0.882}
    \Vertex(394,-144){8.49}
    \Vertex(324,-145){8.49}
    \Vertex(361,-235){9.9}
    \Line(363,-234)(273,-175)
    \Vertex(270,-171){8.49}
  \end{picture}}}\,}

\def\haex{\,{\scalebox{0.3}{
   \begin{picture}(324,355) (0,0)
    \SetWidth{2}
    \SetColor{Black}
    \GOval(210,235)(30,30)(0){0.882}
    \GOval(210,175)(30,30)(0){0.882}
    \Vertex(240,10){9.9}
    \Line(240,10)(240,55)
    \GOval(240,85)(30,30)(0){0.882}
    \Vertex(210,85){8.49}
    \Vertex(270,85){8.49}
    \Line(240,10)(315,85)
    \Vertex(315,85){8.49}
    \Line(240,10)(165,85)
    \Vertex(165,85){8.49}
    \Line(165,85)(120,145)
    \Line(165,85)(165,145)
    \Line(165,85)(210,145)
    \Vertex(165,145){8.49}
    \GOval(120,175)(30,30)(0){0.882}
    \Vertex(90,175){8.49}
    \Vertex(150,175){8.49}
    \Vertex(180,175){8.49}
    \Vertex(210,205){8.49}
    \Vertex(240,175){8.49}
    \Vertex(180,250){8.49}
    \Vertex(240,250){8.49}
  \end{picture}}}\,}

\newcommand{\tree}{\hskip 0.8pc\scalebox{-0.3}{{\parbox{0.5pc}{
  \begin{picture}(30,45) (75,-60)
    \SetWidth{1.5}
    \SetColor{Black}
    \Line(90,-30)(75,-60)
    \Line(90,-30)(105,-60)
    \Line(90,-15)(90,-30)
  \end{picture}}}}}

\newcommand{\treeA}{\hskip 1.5pc\scalebox{-0.3}{{\parbox{0.5pc}{
   \begin{picture}(60,75) (75,-30)
    \SetWidth{1.5}
    \SetColor{Black}
    \Line(90,0)(75,-30)
    \Line(90,0)(105,-30)
    \Line(105,30)(90,0)
    \Line(105,30)(135,-30)
    \Line(105,45)(105,30)
  \end{picture}}}}}

\newcommand{\treeB}{\hskip 1.5pc\scalebox{-0.3}{{\parbox{0.5pc}{
  \begin{picture}(60,75) (75,-30)
    \SetWidth{1.5}
    \SetColor{Black}
    \Line(90,0)(75,-30)
    \Line(105,30)(135,-30)
    \Line(105,45)(105,30)
    \Line(120,0)(105,-30)
    \Line(105,30)(90,0)
  \end{picture}}}}}

\newcommand{\treeC}{\hskip 1.5pc\scalebox{-0.2}{{\parbox{0.5pc}{
 \begin{picture}(90,105) (75,-30)
    \SetWidth{2.4}
    \SetColor{Black}
    \Line(90,0)(75,-30)
    \Line(120,60)(90,0)
    \Line(90,0)(105,-30)
    \Line(105,30)(135,-30)
    \Line(120,60)(165,-30)
    \Line(120,75)(120,60)
  \end{picture}
}}}}

\newcommand{\treeD}{\hskip 1.5pc\scalebox{-0.2}{{\parbox{0.5pc}{
 \begin{picture}(90,105) (75,-30)
    \SetWidth{2.4}
    \SetColor{Black}
    \Line(90,0)(75,-30)
    \Line(120,60)(90,0)
    \Line(90,0)(105,-30)
    \Line(120,60)(165,-30)
    \Line(120,75)(120,60)
    \Line(150,0)(135,-30)
  \end{picture}
}}}}

\newcommand{\treeE}{\hskip 1.5pc\scalebox{-0.2}{{\parbox{0.5pc}{
 \begin{picture}(90,105) (75,-30)
    \SetWidth{2.4}
    \SetColor{Black}
    \Line(90,0)(75,-30)
    \Line(120,60)(90,0)
    \Line(120,60)(165,-30)
    \Line(120,75)(120,60)
    \Line(150,0)(135,-30)
    \Line(135,30)(105,-30)
  \end{picture}
}}}}

\newcommand{\treeF}{\hskip 1.5pc\scalebox{-0.2}{{\parbox{0.5pc}{
 \begin{picture}(90,105) (75,-30)
    \SetWidth{2.4}
    \SetColor{Black}
    \Line(90,0)(75,-30)
    \Line(120,60)(90,0)
    \Line(120,60)(165,-30)
    \Line(120,75)(120,60)
    \Line(135,30)(105,-30)
    \Line(120,0)(135,-30)
  \end{picture}
}}}}

\newcommand{\treeG}{\hskip 1.5pc\scalebox{-0.2}{{\parbox{0.5pc}{
 \begin{picture}(90,105) (75,-30)
    \SetWidth{2.4}
    \SetColor{Black}
    \Line(90,0)(75,-30)
    \Line(120,60)(90,0)
    \Line(120,60)(165,-30)
    \Line(120,75)(120,60)
    \Line(105,30)(135,-30)
    \Line(120,0)(105,-30)
  \end{picture}
}}}}

\newcommand{\treebplus}{\scalebox{0.3}{{\parbox{0.5pc}{
\begin{picture}(457,327) (72,-58)
    \SetWidth{2}
    \SetColor{Black}
    \Line(97,118)(187,41)
    \Line(188,39)(188,-58)
    \Line(188,39)(373,145)
    \DashLine(374,146)(424,173){2}
    \Line(262,82)(200,144)
    \Line(343,128)(290,181)
    \Line(461,197)(423,235)
    \Line(424,174)(529,237)
    \Text(72,139)[lb]{\Huge{\Black{$\Phi^{-1}(t_1)$}}}
    \Text(175,171)[lb]{\Huge{\Black{$\Phi^{-1}(t_2)$}}}
    \Text(272,197)[lb]{\Huge{\Black{$\Phi^{-1}(t_3)$}}}
    \Text(401,253)[lb]{\Huge{\Black{$\Phi^{-1}(t_n)$}}}
  \end{picture}
}}}}

\newcommand{\treebc}{\scalebox{0.21}{{
    \begin{picture}(100,105) (100,-90)
    \SetWidth{1.0}
    \SetColor{Black}
    \Line(150,-135)(150,-150)
    \Line(105,-45)(150,-135)
    \Line(195,-45)(150,-135)
    \Line(165,-45)(165,-105)
    \Line(135,-45)(165,-105)
  \end{picture}
}}}

\newcommand{\treec}{\scalebox{0.21}{{\parbox{0.5pc}{
  \begin{picture}(60,75) (135,-140)
    \SetWidth{1.0}
    \SetColor{Black}
    \Line(165,-75)(165,-120)
    \Line(135,-75)(165,-135)
    \Line(165,-105)(165,-150)
    \Line(195,-75)(165,-135)
  \end{picture}
}}}}

\newcommand{\treebcbis}{\scalebox{0.1}{{\parbox{0.5pc}{
  \begin{picture}(175,309) (78,-91)
    \SetWidth{5}
    \SetColor{Black}
    \Line(78,218)(120,115)
    \Line(120,116)(162,13)
    \Line(162,12)(162,-91)
    \Line(253,216)(162,13)
    \Line(162,218)(162,16)
    \Line(128,218)(120,116)
  \end{picture}
}}}}

\newcommand{\treecbsym}{\scalebox{0.1}{{\parbox{0.5pc}{
 \begin{picture}(175,309) (78,-91)
    \SetWidth{5}
    \SetColor{Black}
    \Line(78,218)(120,115)
    \Line(120,116)(162,13)
    \Line(162,12)(162,-91)
    \Line(253,216)(162,13)
    \Line(187,218)(163,116)
    \Line(140,218)(161,115)
    \Line(162,116)(162,16)
  \end{picture}
}}}}

\newcommand{\treecb}{\scalebox{0.1}{{\parbox{0.5pc}{
\begin{picture}(175,309) (78,-91)
    \SetWidth{5}
    \SetColor{Black}
    \Line(120,218)(120,115)
    \Line(78,218)(120,115)
    \Line(162,218)(120,115)
    \Line(120,116)(162,13)
    \Line(162,12)(162,-91)
    \Line(253,216)(162,13)
  \end{picture}
}}}}

\newcommand{\treecbbis}{\scalebox{0.1}{{\parbox{0.5pc}{
\begin{picture}(175,309) (78,-91)
    \SetWidth{5}
    \SetColor{Black}
    \Line(78,218)(120,115)
    \Line(120,116)(162,13)
    \Line(162,12)(162,-91)
    \Line(253,216)(162,13)
    \Line(162,218)(162,16)
    \Line(199,217)(207,114)
  \end{picture}
}}}}  

 \newcommand{\treed}{\scalebox{0.1}{{\parbox{0.5pc}{ 
\begin{picture}(272,343) (134,-119)
    \SetWidth{5}
    \SetColor{Black}
    \Line(406,222)(270,14)
    \Line(315,224)(270,15)
    \Line(225,222)(270,17)
    \Line(270,-119)(270,14)
    \Line(134,224)(271,15)
  \end{picture}
}}}} 

\newcommand{\treeex}{\scalebox{0.20}{{\parbox{0.5pc}{ 
\begin{picture}(358,475) (91,-103)
    \SetWidth{3}
    \SetColor{Black}
    \Line(91,372)(270,56)
    \Line(119,370)(186,205)
    \Line(185,208)(211,236)
    \Line(150,370)(210,236)
    \Line(211,236)(222,257)
    \Line(179,370)(222,258)
    \Line(239,370)(222,260)
    \Line(209,368)(231,318)
    \Line(269,369)(232,319)
    \Line(300,368)(223,259)
    \Line(420,369)(429,336)
    \Line(449,369)(269,57)
    \Line(389,369)(385,258)
    \Line(269,57)(269,-103)
    \Line(358,370)(384,259)
    \Line(330,369)(222,259)
  \end{picture}
}}}}

\newcommand{\labtreeex}{\scalebox{0.20}{{\parbox{0.5pc}{ 
\begin{picture}(358,475) (91,-103)
    \SetWidth{3}
    \SetColor{Black}
    \Line(91,372)(270,56)
    \Line(119,370)(186,205)
    \Line(185,208)(211,236)
    \Line(150,370)(210,236)
    \Line(211,236)(222,257)
    \Line(179,370)(222,258)
    \Line(239,370)(222,260)
    \Line(209,368)(231,318)
    \Line(269,369)(232,319)
    \Line(300,368)(223,259)
    \Line(420,369)(429,336)
    \Line(449,369)(269,57)
    \Line(389,369)(385,258)
    \Line(269,57)(269,-103)
    \Line(358,370)(384,259)
    \Line(330,369)(222,259)
\Text(447,380)[lb]{\Huge{\Black{$1$}}}
\Text(415,380)[lb]{\Huge{\Black{$2$}}}
\Text(385,380)[lb]{\Huge{\Black{$3$}}}
\Text(353,380)[lb]{\Huge{\Black{$4$}}}
\Text(325,380)[lb]{\Huge{\Black{$5$}}}
\Text(295,380)[lb]{\Huge{\Black{$6$}}}
\Text(264,380)[lb]{\Huge{\Black{$7$}}}
\Text(234,380)[lb]{\Huge{\Black{$8$}}}
\Text(204,380)[lb]{\Huge{\Black{$9$}}}
\Text(170,380)[lb]{\Huge{\Black{$10$}}}
\Text(140,380)[lb]{\Huge{\Black{$11$}}}
\Text(109,380)[lb]{\Huge{\Black{$12$}}}
\Text(77,380)[lb]{\Huge{\Black{$13$}}}

  \end{picture}
}}}}

\def\labhaex{\,{\scalebox{0.3}{
   \begin{picture}(324,355) (0,0)
    \SetWidth{2}
    \SetColor{Black}
    \GOval(210,235)(30,30)(0){0.882}
    \GOval(210,175)(30,30)(0){0.882}
    \Vertex(240,10){9.9}
    \Line(240,10)(240,55)
    \GOval(240,85)(30,30)(0){0.882}
    \Vertex(210,85){8.49}
    \Vertex(270,85){8.49}
    \Line(240,10)(315,85)
    \Vertex(315,85){8.49}
    \Line(240,10)(165,85)
    \Vertex(165,85){8.49}
    \Line(165,85)(120,145)
    \Line(165,85)(165,145)
    \Line(165,85)(210,145)
    \Vertex(165,145){8.49}
    \GOval(120,175)(30,30)(0){0.882}
    \Vertex(90,175){8.49}
    \Vertex(150,175){8.49}
    \Vertex(180,175){8.49}
    \Vertex(210,205){8.49}
    \Vertex(240,175){8.49}
    \Vertex(180,250){8.49}
    \Vertex(240,250){8.49}
\Text(255,0)[lb]{\huge{\Black{$1$}}}
\Text(335,85)[lb]{\huge{\Black{$2$}}}
\Text(285,85)[lb]{\huge{\Black{$3$}}}
\Text(188,84)[lb]{\huge{\Black{$4$}}}
\Text(140,80)[lb]{\huge{\Black{$5$}}}
\Text(260,175)[lb]{\huge{\Black{$6$}}}
\Text(240,200)[lb]{\huge{\Black{$7$}}}
\Text(260,250)[lb]{\huge{\Black{$8$}}}
\Text(150,250)[lb]{\huge{\Black{$9$}}}
\Text(165,192)[lb]{\huge{\Black{$10$}}}
\Text(176,132)[lb]{\huge{\Black{$11$}}}
\Text(141,199)[lb]{\huge{\Black{$12$}}}
\Text(52,175)[lb]{\huge{\Black{$13$}}}
  \end{picture}}}\,}

\def\pentagone{\,{\scalebox{0.6}{
   \begin{picture}(512,380) (76,-90)
    \SetWidth{1}
    \SetColor{Black}
    \Line(225,230)(455,230)
    \Line(455,230)(530,76)
    \Line(530,76)(343,-56)
    \Line(343,-56)(155,76)
    \Line(155,76)(224,230)
    \Text(190,261)[lb]{\Large{\Black{$\treeE$}}}
    \Text(448,274)[lb]{\Large{\Black{$\treeF$}}}
    \Text(322,-90)[lb]{\Large{\Black{$\treeC$}}}
    \Text(76,78)[lb]{\Large{\Black{$\treeD$}}}
    \Text(558,84)[lb]{\Large{\Black{$\treeG$}}}
    \Text(311,104)[lb]{\Large{\Black{$\treed$}}}
    \Text(319,266)[lb]{\Large{\Black{$\treecb$}}}
    \Text(521,161)[lb]{\Large{\Black{$\treecbsym$}}}
    \Text(461,-18)[lb]{\Large{\Black{$\treebc$}}}
    \Text(109,174)[lb]{\Large{\Black{$\treebcbis$}}}
    \Text(186,-45)[lb]{\Large{\Black{$\treecbbis$}}}
  \end{picture}}}\,}
  
  \def\pentagonebis{\,{\scalebox{0.6}{
   \begin{picture}(512,380) (76,-90)
    \SetWidth{1}
    \SetColor{Black}
    \Line(225,230)(455,230)
    \Line(455,230)(530,76)
    \Line(530,76)(343,-56)
    \Line(343,-56)(155,76)
    \Line(155,76)(224,230)
    \Text(190,261)[lb]{\Black{$\arbrecaB$}}
    \Text(448,274)[lb]{\Black{$\arbrecbB$}}
    \Text(322,-90)[lb]{\Black{$\arbrecdB$}}
    \Text(76,78)[lb]{\Black{$\arbreccB$}}
    \Text(558,84)[lb]{\Black{$\arbreceB$}}
    \Text(311,104)[lb]{\Black{$\hac$}}
    \Text(319,266)[lb]{\Black{$\hacab$}}
    \Text(521,161)[lb]{\Black{$\hacbe$}}
    \Text(461,-18)[lb]{\Black{$\haced$}}
    \Text(109,174)[lb]{\Black{$\hacca$}}
    \Text(186,-45)[lb]{\Black{$\hacdc$}}
  \end{picture}}}\,}

\nc{\ignore}[1]{{}}
\nc{\mrm}[1]{{\rm #1}}
\nc{\dirlim}{\displaystyle{\lim_{\longrightarrow}}\,}
\nc{\invlim}{\displaystyle{\lim_{\longleftarrow}}\,}
\nc{\vep}{\varepsilon} \nc{\ep}{\epsilon}
\nc{\sigmat}{\widetilde\sigma}
\nc{\ostar}{\overline{*}}

\nc{\mchar}{\mrm{Char}}
\nc{\Hom}{\mrm{Hom}}
\nc{\id}{\mrm{id}}

\nc{\remark}{\noindent{\bf{Remark:}}}
\nc{\remarks}{\noindent{\bf{Remarks:}}}

 \nc{\delete}[1]{}
 \nc{\grad}[1]{^{({#1})}}
 \nc{\fil}[1]{_{#1}}

\nc{\BA}{{\Bbb A}} \nc{\CC}{{\Bbb C}} \nc{\DD}{{\Bbb D}}
\nc{\EE}{{\Bbb E}} \nc{\FF}{{\Bbb F}} \nc{\GG}{{\Bbb G}}
\nc{\HH}{{\Bbb H}} \nc{\LL}{{\Bbb L}} \nc{\NN}{{\Bbb N}}
\nc{\PP}{{\Bbb P}} \nc{\QQ}{{\Bbb Q}} \nc{\RR}{{\Bbb R}}
\nc{\TT}{{\Bbb T}} \nc{\VV}{{\Bbb V}} \nc{\ZZ}{{\Bbb Z}}
\nc{\Cal}[1]{{\mathcal {#1}}}
\nc{\mop}[1]{\mathop{\hbox {\rm #1} }\nolimits}
\nc{\smop}[1]{\mathop{\hbox {\eightrm #1} }\nolimits}
\nc{\mopl}[1]{\mathop{\hbox {\rm #1} }\limits}
\nc{\frakg}{{\frak g}}
\nc{\g}[1]{{\frak {#1}}}
\def \restr#1{\mathstrut_{\textstyle |}\raise-8pt\hbox{$\scriptstyle #1$}}
\def \srestr#1{\mathstrut_{\scriptstyle |}\hbox to
  -1.5pt{}\raise-4pt\hbox{$\scriptscriptstyle #1$}}
\nc{\wt}{\widetilde}
\nc{\wh}{\widehat}
\nc{\un}{\hbox{\bf 1}}
\nc{\redtext}[1]{\textcolor{red}{\tt #1}}
\nc{\bluetext}[1]{\textcolor{blue}{#1}}
\nc{\comment}[1]{[[{\tt {#1}}]] }
\nc{\R}{\mathbb R}
\nc\fleche[1]{\mathop{\hbox to #1 mm{\rightarrowfill}}\limits}
\def\semi{\mathrel{\times}\kern -.85pt\joinrel\mathrel{\raise
    1.4pt\hbox{${\scriptscriptstyle |}$}}}
\nc{\np}{/\hskip -2.3mm\pi}
\nc{\snp}{/\hskip -1.8mm\pi}

\def\ta1{{\scalebox{0.2}{ %%%%%%%%%%%%%%%%%%%%%%%%%%%%%%%%%\ta1
\begin{picture}(12,12)(38,-38)
\SetWidth{0.5} \SetColor{Black} \Vertex(45,-33){5.66}
\end{picture}}}}
\begin{document}

\title[Rooted trees and hypertrees]
      {On an extension of Knuth's rotation \\correspondence to reduced planar trees}

\author{Kurusch Ebrahimi-Fard}
\address{Instituto de Ciencias Matem\'aticas,
		C/ Nicol\'as Cabrera, no.~13-15, 28049 Madrid, Spain.
		On leave from Univ.~de Haute Alsace, Mulhouse, France}
         \email{kurusch@icmat.es, kurusch.ebrahimi-fard@uha.fr}         
         \urladdr{www.icmat.es/kurusch}

\author{Dominique Manchon}
\address{Univ.~Blaise Pascal,
         	C.N.R.S.-UMR 6620,
         	63177 Aubi\`ere, France}       
         \email{manchon@math.univ-bpclermont.fr}
         \urladdr{http://math.univ-bpclermont.fr/~manchon/}

%%%%%%%%%%%%%%%%%%%%%%%%%%%%%%%%%%%%%%%%%%%%%%%%%%%%%%%%%%%%%%%%%%%
\date{March 2nd, 2012}
%%%%%%%%%%%%%%%%%%%%%%%%%%%%%%%%%%%%%%%%%%%%%%%%%%%%%%%%%%%%%%%%%%%

\begin{abstract}
We present a bijection from planar reduced trees to planar rooted hypertrees, which extends Knuth's rotation correspondence between planar binary trees and planar rooted trees. The operadic counterpart of the new bijection is explained. Related to this, the space of planar reduced forests is endowed with a combinatorial Hopf algebra structure. The corresponding structure on the space of planar rooted hyperforests is also described.
\end{abstract}

%%%%%%%%%%%%%%%%%%%%%%%%%%%%%%%%%%%%%%%%%%%%%%%%%%%%%%%%%%%%%%%%%%%

\maketitle

%%%%%%%%%%%%%%%%%%%%%%%%%%%%%%%%%%%%%%%%%%%%%%%%%%%%%%%%%%%%%%%%%%%

\tableofcontents

%%%%%%%%%%%%%%%%%%%%%%%%%%%%%%%%%%%%%%%%%%%%%%%%%%%%%%%%%%%%%%%%%%%

\section{Introduction}
\label{sect:intro}

Rooted trees have been extensively used in many branches of pure and applied mathematics. Especially in the latter case they gained particular prominence due to the pioneering work on numerical integration methods by John Butcher in the 1960s \cite{Butcher1,Duer,HW}. He discovered a group structure in the context of Runge--Kutta integration methods. This group structure encodes the composition of so-called $B$-series. The latter are a generalization of Taylor series, in which rooted trees naturally appear, as Arthur Cayley noticed in his classical 1857 paper \cite{Cayley}. See \cite{Butcher2,HWL} for details. Since then, algebraic structures have become an important aspect in the study of numerical methods and related fields, see e.g.~\cite{BeOw,CMOQ,HNT,Iserles,Murua}. 

Somewhat after Butcher's seminal work, Gian-Carlo Rota and Saj-Nicole Joni observed in a seminal paper \cite{J-R}, that various combinatorial objects naturally possess compatible product and coproduct structures. With the work by William Schmitt \cite{Schmitt} this ultimately converged into the notion of combinatorial Hopf algebra, i.e., as Marcelo Aguiar puts it, a connected graded vector space where the homogeneous components are spanned by finite sets of combinatorial objects, and the algebraic structures are given by particular constructions on those objects. Rooted trees provide a genuine example for such combinatorial objects, and several Hopf algebra structures have been described using them. In \cite{GL, Hoffman,Holtkamp1,Holtkamp2,Zhao} the reader finds more details. In particular, Arne D\"ur, and later Christian Brouder \cite{Brouder} showed that the Butcher group identifies with the group of characters on the dual of a commutative graded Hopf algebra of rooted trees described by Alain Connes and Dirk Kreimer \cite{CK}. In \cite{CEM} a combinatorial Hopf algebra structure on rooted trees was described that corresponds to the substitution law of B-series introduced in \cite{CHV1}, see also \cite{CHV2}. In \cite{MKW} a Hopf algebra on planar rooted trees was introduced in the context of Lie--Butcher series on Lie groups. 

Combinatorial Hopf algebras on rooted trees are generally related to the fact that free pre-Lie algebras are naturally described in terms of rooted trees \cite{ChaLiv,DL,Segal}. In the case of Hans Munthe-Kaas' and William Wright's noncommutative Hopf algebra for Lie--Butcher series \cite{MKW}, this has been generalized to so-called $D$-algebras. Fr\'ed\'eric Chapoton observed in \cite{Chap1} that an operadic approach may provide an adequate perspective on the link between pre-Lie structures, the group of characters and combinatorial Hopf algebras.

The theory of correspondences between combinatorial objects is one of the main topics in combinatorics. As an example we mention Robinson's and Schensted's bijection between permutations and standard tableaux. Another example is Donald Knuth's rotation correspondence \cite{Knuth68} for planar binary trees, which maps a planar binary tree with $n-1$ internal vertices into a planar rooted tree with $n$ vertices. In this paper we generalize Knuth's correspondence to a bijection between planar reduced trees and planar rooted hypertrees. This bijection is used to transfer a combinatorial Hopf algebra structure on planar reduced trees to planar rooted hypertrees. It turns out that the coproduct of the latter is very similar to the one in Munthe-Kaas and Wright's Hopf algebra. In a forthcoming article we will describe in more detail the underlying reason for this.  

\medskip 

This paper is organized as follows. In Section \ref{sect:planar} we introduce the notions of planar binary trees and planar rooted trees. The Butcher product on trees is presented. We recall Knuth's rotation correspondence between planar binary trees and planar rooted trees, which we then extend to a bijection from planar reduced trees to so-called planar rooted hypertrees. In Section \ref{sect:Hopf} we briefly recall some notions from the theory of connected graded bialgebras, and then define a Hopf algebra on planar reduced trees respectively planar rooted hypertrees.

%%%%%%%%%%%%%%%%%%%%%%%%%%%%%%%%%%%%%%%%%%%%%

\section{Planar rooted trees and hypertrees}
\label{sect:planar}

Recall that a tree is an undirected connected graph made out of vertices and edges. It is without cycles, that is, any two vertices can be connected by exactly one simple path. We denote the set of vertices and edges of a tree $t$ by $V(t)$ and $E(t)$, respectively. In this section we introduce the objects of this work, which are particular classes of trees, i.e.~planar binary (reduced) trees and planar rooted (hyper)trees.

%%%%%%%%%%%%%%%%%%%%%%%%%%%%%%%%%%%%%%%%%%%%%

\subsection{Planar trees}
\label{ssect:pbt}

We start with the notion of a {\it{planar binary}} tree, which is a finite oriented tree given an embedding in the plane, such that all vertices have exactly two incoming edges and one outgoing edge. An edge can be internal (connecting two vertices) or external (with one loose end). The external incoming edges are the leaves. The root is the unique edge not ending in a vertex. 
$$
	\vert
	\quad\
	\scalebox{1.5}{\tree}
	\quad
	\treeA
	\quad\
	\treeB
	\quad\
	\treeE
	\quad
	\treeD
	\quad
	\treeC
	\quad
	\treeG	
	\quad 
	\treeF \quad \ldots
$$
The single edge $\vert$ is the unique planar binary tree without internal vertices. We denote by $T^{bin}_{pl}$ (resp.~$\Cal T^{bin}_{pl}$) the set (resp.~the linear span) of planar binary trees. A simple grading for such trees is given in terms of the number of internal vertices. Alternatively, one can use the number of leaves. Observe that for any pair of planar binary trees $t_1,t_2$ we can build up a new planar binary tree via the grafting operation, $t_3:=t_1 \vee t_2$, i.e.~by considering the $Y$-shaped tree $\hspace{-0.2cm}\begin{array}{cc}&\\[-0.7cm]\scalebox{0.8}{\tree}\hspace{-0.5cm}\end{array}$ (the unique planar binary tree with two leaves) and replacing the left branch (resp.~the right branch) by $t_1$ (resp.~$t_2$). 
$$
	\vert \vee \vert = \tree
	\quad\
	\tree \vee \vert = \treeB
	\quad\
	\vert \vee \tree = \treeA
	\quad\
	\tree \vee \tree = \treeD
	\quad\
	\vert \vee \treeB = \treeG
$$

Seen as a product on $\Cal T^{bin}_{pl}$, the grafting operation $\vee$ is neither associative nor commutative, $t_1 \vee t_2 \neq t_2 \vee t_1$. In fact, one can show that it is purely magmatic. Notice that this product is of degree one with respect to the grading in terms of internal vertices, i.e.~ for two trees $t_1,t_2$ of degrees $n_1,n_2$, respectively, the product $t_1 \vee t_2$ is of degree $n_1+n_2+1$. However, with respect to the leave number grading this product is of degree zero. 

A {\it{planar rooted tree}} is a finite oriented rooted tree given an embedding in the plane, such that all vertices, except one, have arbitrarily many incoming edges and one outgoing edge. The root is the one vertex without an outgoing edge. 
$$
	\racine
	\quad\
	\arbrea
	\quad
	\arbrebb
	\quad\
	\arbreba
	\quad\
	\arbreca
	\quad
	\arbrecc
	\quad
	\arbrecd
	\quad\
	\arbrece	
	\quad
	\arbrecb \quad\ \cdots
$$

The single vertex $\racine$ is the unique rooted tree without edges. Note that we put the root at the bottom of the tree. The set (resp.~the linear span) of planar non-empty rooted trees will be denoted by $T_{pl}$ (resp.~$\Cal T_{pl}$). A natural grading for such trees is given in terms of the number of edges. Another one is given by the number of vertices. Observe that any rooted tree of degree bigger than zero writes in a unique way:
\begin{equation*}
	t=B_+(t_1\cdots t_n),
\end{equation*}
where $B_+$ associates to the forest $t_1\cdots t_n$ the planar tree obtained by grafting all the planar trees $t_j$, $j=1,\ldots,n$, on a common root. 
$$
	B_+(\racine)=\arbrea
	\quad\
	B_+(\racine\racine)=	\arbrebb
	\quad\
	B_+(\arbrea\racine)=\arbrecc
	\quad\
	B_+(\racine\arbrea)=\arbrece.
$$
Recall that sometimes, one finds the notation $t=[t_1\cdots t_n]$ in the literature \cite{Brouder,Butcher1}. Note that the order in which the branch trees are displayed has to be taken into account. 

Further below we will recall the classical correspondence between these two types of trees, due to Knuth \cite{Knuth68}.

%%%%%%%%%%%%%%%%%%%%%%%%%%%%%%%%%%%%%%%%%%%%%

\subsection{The Butcher product}
\label{ssect:prtbutcher}

Motivated by the use of (non-)planar rooted trees in the theory of numerical integration methods \cite{Butcher1,Butcher2,GL,HW}, we introduce a planar version of the classical Butcher product. The (left) {\sl{Butcher product}} of two planar rooted trees $t,u$ is defined by connecting the root of $t$ via a new edge to the root of $u$ such that $t$  becomes the leftmost branch tree, that is, for two trees $t=B_+(t_1\cdots t_n)$ and $u=B_+(u_1\cdots u_p)$:
\begin{equation}
\label{butcherProd}
	t \butcher u := B_+(t u_1\cdots u_p).
\end{equation} 
Observe that it is neither associative nor commutative, and, again contrarily to the non-planar case, it is also not NAP (Non-Associative Permutative) \cite{Li1}, i.e.~it does not satisfy the identity $t \butcher (u \butcher v) = u \butcher (t \butcher v)$.

%%%%%%%%%%%%%%%%%%%%%%%%%%%%%%%%%%%%%%%%%%%%%

\subsection{Knuth's correspondence between planar binary and planar rooted trees}
\label{ssect:pbt2prt}

Knuth describes in \cite{Knuth68} a natural relation between planar rooted trees and planar binary trees, known as rotation correspondence. We only give a recursive definition of this bijection, and refrain from providing more details. The interested reader is refered to Marckert's paper \cite{Marckert} for a nice description of the rotational aspect. 

Recall that by the single edge $\vert$ we denote the unique planar binary tree without internal vertices. Now we recursively define a map $\Phi:T^{bin}_{pl} \to T_{pl}$ by $\Phi(\vert):=\racine$ and:
\begin{equation}
\label{phi}
	\Phi(t_1\vee t_2):=\Phi(t_1)\butcher\Phi(t_2).
\end{equation}
This map is clearly well-defined and bijective\footnote{It appears also in \cite{Foissy} in slightly different form.}, with its inverse recursively given by:

\begin{equation}
\label{phi-inv}
	\hspace{-2cm}\Phi^{-1}\big(B_+(t_1\cdots t_n)\big) = {\scalebox{0.8}{\treebplus}}
\end{equation}
The first few terms write:
$$
	\Phi(|)=\racine
	\qquad\
	\Phi(\!\tree)=\arbrea
	\qquad
	\Phi({\scalebox{0.8}{\treeA}})=\arbrebb
	\qquad\
	\Phi({\scalebox{0.8}{\treeB}})=\arbreba
$$
$$
	\Phi({\scalebox{0.8}{\treeE}})=\arbreca
	\qquad\
	\Phi({\scalebox{0.8}{\treeD}})=\arbrecc
	\qquad\
	\Phi({\scalebox{0.8}{\treeC}})=\arbrecd
	\qquad\
	\Phi({\scalebox{0.8}{\treeG}})=\arbrece	
	\qquad\ 
	\Phi({\scalebox{0.8}{\treeF}})=\arbrecb.
$$
\vspace{0.2cm}

Note that this simple bijection implies that the left Butcher product (\ref{butcherProd}) is also magmatic.

%%%%%%%%%%%%%%%%%%%%%%%%%%%%%%%%%%%%%%%%%%%%%

\subsection{Reduced planar rooted trees and planar rooted hypertrees}
\label{ssect:hypertrees}

A planar tree is called {\sl{reduced}}, if any inner vertex has at least two incoming edges. We denote by $T^{red}_{pl}$ (resp.~$\Cal T^{red}_{pl}$) the set (resp.~the linear span) of reduced planar trees. Any reduced planar tree can be described for $n>1$ as $t=\bigvee(t_1,\ldots,t_n)$, i.e.~it can be obtained by considering the unique tree with one internal vertex and $n$ incoming edges, and replacing the $i^{\smop{th}}$ branch by $t_i$. There is a partial order on $T^{red}_{pl}$ defined as follows: $t_1 \le t_2$ if $t_1$ can be obtained from $t_2$ by glueing some inner vertices together. In particular, two comparable trees must have the same number of leaves. The minimal elements are the trees with only one inner vertex, and the maximal elements are the planar binary trees.

\medskip

We would like to propose a way to extend the bijection $\Phi$, originally defined on planar binary trees, to reduced planar trees, thus answering a question by J.-L.~Loday. The image of $T^{red}_{pl}$ will be the space $HT_{pl}$ of {\sl planar rooted hypertrees\/}, which we introduce now.

\smallskip

Following Chapoton \cite{Chap2}, a {\sl hypergraph\/} on a finite set $I$ of vertices is a nonempty set of parts of $I$ of cardinality at least $2$, which will be called the {\sl edges\/} of the hypergraph. A {\sl path\/} in a hypergraph is a sequence  $i_1,\ldots ,i_k$ of vertices such that any pair $\{i_j , i_{j+1}\}$ is included in an edge. A hypergraph is connected if any two vertices can be joined by a path. A {\sl hypertree \/} is a connected hypergraph without cycles except those which are included in a single edge. Two different edges in a hypertree then meet at one single vertex or have empty intersection.

\smallskip

A {\sl rooted hypertree\/} is a hypertree with a distinguished vertex. This defines a partial order on the set of edges as follows: $e < e'$ if for any vertex $j$ in $e$ and any vertex  $j'$ in $e'$ there is a path from the root to $j'$ through $j$. This in turn defines a preorder on the vertices in an obvious way. For any edge $e'$ not containing the root, there is a unique edge $e$ such that $e < e'$ and $e \cap e' \not =\emptyset$. The unique vertex in this intersection will be called the root of the edge $e'$.  Define a {\sl planar rooted hypertree\/} as a rooted hypertree together with an embedding into the plane such that any edge is embedded in the boundary of a small topological disk. This defines a partial order on the vertices compatible with the preorder defined above, i.e. it determines a total order on each edge with the edge's root as minimal element, by running counterclockwise along the boundary of the disk. The following planar rooted hypertree:
\vspace{-0.5cm}
\begin{equation}
\label{ex:example}
	{\scalebox{0.8}{\haex}}
\end{equation}
has seven edges altogether, three with 2 vertices, three with 3 vertices and one with 4 vertices. Each edge of cardinality bigger than 2 is represented by a blob. The vertices are drawn on the circle delimiting the blob, and are ordered counterclockwise starting from the edge's root.

\smallskip

There is a partial order on the set of all rooted planar hypertrees on a given set $I$ of vertices with root $r\in I$ fixed: $t_1\le t_2$ if and only if any edge of $t_2$ is contained in an edge of $t_1$.  The minimal element is the hypertree with only one edge equal to the whole $I$, and the maximal elements are planar rooted trees on $I$ with root $r$.

\medskip

We are now ready to extend the bijection $\Phi$. For any ordered collection $(t_1,\ldots, t_n)$ of planar rooted hypertrees with respective roots $r_j$ we define $\beta(t_1,\ldots t_n)$ by collecting the roots $r_j$ into a common edge, in which the vertices are put in the {\sl reversed order\/}. In particular, it implies that $r_n$ is the root of this new edge, hence the root of the new built hypertree. This certainly extends the Butcher product of two trees (\ref{butcherProd}). We then extend $\Phi$ by setting recursively:
\begin{equation}
\label{eq:phi}
	\Phi\big(\bigvee(t_1,\ldots,t_n)\big):=\beta\big(\Phi(t_1),\ldots,\Phi(t_n)\big).
\end{equation}
Any planar rooted hypertree writes in a unique way as $\beta(s_1,\ldots,s_n)$, where $n$ is the cardinality of the leftmost edge containing the root. The inverse $\Phi^{-1}$ is then recursively defined as follows:
\begin{equation}
\label{phi-invbis}
	\Phi^{-1}\big(\beta(s_1,\ldots,s_n)\big)=\bigvee\big(\Phi^{-1}(s_1),\ldots,\Phi^{-1}(s_n)\big).
\end{equation}
Considering the example (\ref{ex:example}) above, we have:
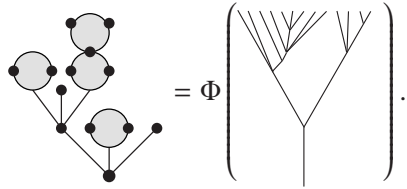
\vspace{-0.5cm}
\begin{equation*}
\raise -11mm\hbox{{${\scalebox{0.8}{\haex}}$}}=\Phi\left({\scalebox{0.7}{\treeex}}\hskip 19mm\right).
\end{equation*}

%%%%%%%%%%%%%%%%%%%%%%%%%%%%%%%%%%%

\ignore{Planar rooted hypertrees write in a unique way:
\begin{equation*}
	t=L(t_1\cdots t_n)
\end{equation*}
where $L$ associates to the ordered hyperforest $t_1\cdots t_n$ the planar
hypertree obtained by grafting all the branches $t_j$, $j=1,\ldots,n$, on a
common root $r$, such that the roots $r_j$ of the branches $t_j$ together with
$r$ build up a new edge of cardinality $n+1$. Note that the order in which the
branches are displayed has to be taken into account due to planarity}
%%%%%%%%%%%%%%%%%%%%%%%%%%%%%%%%%%%

Recall that the reduced planar trees with $n$ leaves are in bijection with the cells of the $n-2$-associahedron. In particular, reduced planar trees with four leaves can be displayed on the pentagon like this:
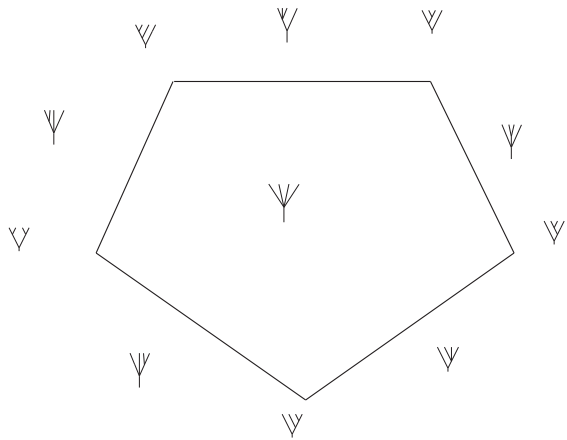
\begin{equation*}
	{\scalebox{0.7}{\pentagone}}
\end{equation*}
Under transformation $\Phi$ the picture transforms like this:
\vskip 3mm
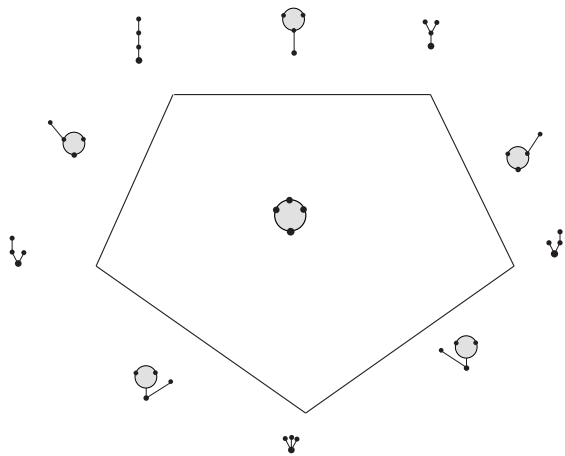
\begin{equation*}
	{\scalebox{0.7}{\pentagonebis}}
\end{equation*}
It is easy to show that $\Phi$ respects the partial orders defined above, which are two manifestations of the reverse incidence order of the associahedron.

%%%%%%%%%%%%%%%%%%%%%%%%%%%%%%%%%%%

\subsection{Adding decorations}
\label{sect:deco}

A planar binary tree decorated by a set $I$ is a planar binary tree together with a map $\delta$ form the set of its internal vertices to $I$. There are grafting operations $\vee_i,\,i\in I$ defined as in the undecorated case, except that the new internal vertex is decorated by
$i$. This decoration procedure generalizes to planar reduced trees as follows: given a partitioned set $\Cal I=I_2\amalg I_3\amalg\cdots$, a planar reduced tree decorated by $\Cal I$ is a planar reduced tree together with a map $\delta$ from the set of its internal vertices to $\Cal I$, which sends the internal vertices which have $n$ incoming edges into $I_n$. Any such decorated planar reduced rooted tree can be uniquely written as:
\begin{equation*}
	t=\bigvee_i(t_1,\ldots,t_n)
\end{equation*}
with $i\in I_n$ for some $n\ge 2$, i.e.~it can be obtained by considering the unique tree with one internal vertex decorated by $i$ and $n$ incoming edges, and replacing the $n^{\smop{th}}$ branch by $t_n$.

\medskip

Equations \eqref{eq:phi} and \eqref{phi-invbis} also recursively define a bijection $\omega_t$ between the internal vertices of a reduced planar tree $t$ and the edges of the planar rooted hypertree $\Phi(t)$, which associates to any internal vertex of $t$ with $n$
incoming edges an edge of $\Phi(t)$ with $n$ vertices. The vertex of $t=\bigvee(t_1,\ldots t_n)$ closest to the root (with $n$ incoming edges) is sent to the leftmost edge of $\Phi(t)$ containing the root. The bijection $\Phi$ hence gives rise to a bijection $\Phi_{\Cal I}$ from $\Cal I$-decorated reduced planar trees to rooted planar hypertrees with edges decorated by $\Cal I$ (i.e.~the edges with $n$ vertices are decorated by $I_n,\,n\ge 2$). The bijection $\Phi_{\Cal I}$ is defined as follows:
\begin{equation*}
	\Phi_{\Cal I}(t,\,\delta)=(\Phi(t),\,\delta\circ\omega_t^{-1}).
\end{equation*}
Any such $\Cal I$-decorated rooted planar hypertree can be uniquely written as:
\begin{equation*}
	s=\beta_i(s_1,\ldots,s_n)
\end{equation*}
with $i\in I_n$ for some $n\ge 2$, i.e.~it can be obtained by collecting the roots of $s_j,\,j=1,\ldots ,n$, this making the leftmost bottom edge, and decorating this new edge by $i$.

%%%%%%%%%%%%%%%%%%%%%%%%%%%%%%%%%%%

\subsection{Operadic structure}
\label{ssect:operadic}

Equation \eqref{phi-invbis} recursively defines a bijection between the vertices of a planar rooted hypertree $s$ and the leaves of the reduced planar tree $\Phi^{-1}(s)$ (the root corresponding to the rightmost leaf). Any labeling of the vertices of $s$ thus corresponds to a labeling of the leaves of $\Phi^{-1}(s)$. On the example above this reads:
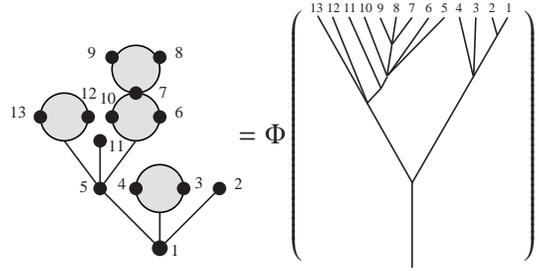
\begin{equation*}
	\raise -15mm\hbox{{$\labhaex$}}=\Phi\left(\hskip 2mm\labtreeex\hskip 27mm\right).
\end{equation*}

Recall that an $\mathbb S$-object is a graded vector space $V=V_1\oplus V_2\oplus\cdots$ together with an action of the symmetric group $S_k$ on $V_k$ for any $k\ge 1$. For any partitioned set $\Cal I=I_2\amalg I_3\amalg\cdots$ we consider the $\mathbb S$-object $V(\Cal I)$ defined by  $(V_{\Cal I})_1=k$ and $(V_{\Cal I})_n=k^{|I_n|}\otimes k[S_n]$ for $n\ge 2$. The vector space $\Cal T^{red,\Cal I}_{pl}$ generated by $\Cal I$-decorated reduced planar trees (see Paragraph \ref{sect:deco} above) naturally encodes the free operad on $V(\Cal I)$: the partial composition $\sigma\circ_i\tau$ of two $\Cal I$-decorated planar reduced trees with labeled leaves is obtained by replacing leaf number $i$ of $\sigma$ by $\tau$. The operadic structure of $\Cal I$-decorated reduced  planar rooted trees can be transferred to the linear span $H\Cal T_{pl}$ of planar rooted hypertrees by means of $\Phi_{\Cal I}^{-1}$: the partial composition $t_1\circ_i t_2$ of two planar rooted hypertrees with labeled vertices is then obtained by replacing vertex number $i$ of hypertree $t_1$ by the root of the hypertree $t_2$, and putting the hypertree $t_2$ plugged this way on the right of the other edges stemming from vertex number $i$. This is easily seen when the vertex is the root of $t_1$, and the other vertices are treated by induction, by remarking that a vertex of $t_1$ different from the root is a vertex of the $j^{\smop{th}}$ branch tree $t_{1,j}$ of $t_1$.

\smallskip

The fully transferred operadic structure on the vector space $H\Cal T_{pl}^{\Cal I}$ of $\Cal I$-decorated hypertrees is then the following: $\gamma(t;t_1,\ldots ,t_n)$ is given by replacing vertex number $i$ of $t$ by the root of $t_i$, and by putting the plugged hypertree $t_i$ on the right. This class of operads is known as generic magmatic operads \cite{Zinbiel}.

%%%%%%%%%%%%%%%%%%%%%%%%%%%%%%%%%%%%%%%%%%%%%

\section{Hopf algebra structures on trees}
\label{sect:Hopf}

%%%%%%%%%%%%%%%%%%%%%%%%%%%%%%%%%%%%%%%%%%%%%

\subsection{Connected filtered bialgebras}
\label{ssect:bialg}

In general, $k$ denotes the ground field (of characteristic zero) over which all algebraic structures are defined. Recall the definition of a  {\it{bialgebra}}, which is an algebra and coalgebra structure together with compatibility relations \cite{J-R}. We denote a {\it{Hopf algebra}} by $(\Cal H,m_{\Cal H},\eta_{\Cal H},\Delta_{\Cal H},\epsilon_{\Cal H},S)$. It is a bialgebra together with a particular $k$-linear map, i.e.~the {\it{antipode}} $S: \Cal H \to \Cal H$, satisfying the Hopf algebra axioms~\cite{Sw69,Manchon}. In the following we omit subscripts if there is no danger of confusion. We denote the unit by $\un=\eta_{\Cal H}(1)$. Let $\Cal H$ be a connected filtered bialgebra, that is:
$$
    k=\Cal H^{(0)} \subset \Cal H^{(1)} \subset \cdots \subset 
      \Cal H^{(n)} \subset \cdots, \hskip 6mm \bigcup_{n\ge 0} \Cal H^{(n)}=\Cal H.
$$
For any $x \in \Cal H^{(n)}$ we have, using a variant of Sweedler's notation \cite{Sw69}:
\begin{equation*}
    \Delta(x) = x \otimes \un + \un \otimes x + \sum_{(x)} x' \otimes x'',
\end{equation*}
where the filtration degrees of $x'$ and $x''$ are strictly smaller than $n$. Recall that by definition we call an element $x
\in \Cal H$ {\it{primitive}} if:
$$
    \bar{\Delta}(x) := \Delta(x) -  x \otimes \un - \un \otimes x = 0.
$$
The antipode $S: \Cal H \to \Cal H$ is defined in terms of the equations:
\begin{equation} 
\label{antipode}
    S  \ast Id =  m \circ (S \otimes Id) \circ \Delta = \eta \circ \epsilon = Id \ast S,
\end{equation}
where the convolution product for two linear maps $ f,g \in \Cal L(\Cal H,\Cal H)$ is defined by $f \ast g := m \circ (f \otimes g) \circ \Delta : \Cal H \to \Cal H$, i.e.: 
\begin{equation*}
    (f \ast g)(x) = f(x) g(\un) + f(\un)g(x) + \sum_{(x)} f(x') g(x'') \in \Cal H.
\end{equation*}
It yields an associative algebra with unit $e:=\epsilon$ on the vector space $\Cal L(\Cal H,\Cal H)$. The antipode always exists for connected filtered bialgebras, hence any connected filtered bialgebra is a {\sl{connected filtered Hopf algebra}\/}. Equations~(\ref{antipode}) imply the following recursive formulas for the antipode starting with $S(\un)=\un$ and for $x \in \mop{ker}\epsilon$:
 \allowdisplaybreaks{
\begin{eqnarray*}
    S(x) = -x- \displaystyle\sum_{(x)}S(x')x'' \ = \ -x- \displaystyle\sum_{(x)}x'S(x'').
\end{eqnarray*}}
Let $\Cal H$ be a graded Hopf algebra. The grading induces a biderivation $ Y: \Cal H^{(n)} \rightarrow \Cal H^{(n)}$ defined on homogeneous elements by $x  \longmapsto nx$.

%%%%%%%%%%%%%%%%%%%%%%%%%%%%%%%%%%%%%%%%%%%%%

\subsection{The Butcher--Connes--Kreimer Hopf algebra of rooted forests}
\label{ssect:CKHopf}

The paradigm of a connected filtered, in fact, graded, Hopf algebra is given by the Butcher--Connes--Kreimer Hopf algebra ${\mathcal H}_{\makebox{{\tiny{BCK}}}}$ of rooted forests over $k$, graded by the number of vertices \cite{Butcher1,CK,Duer,Manchon}. It is the free unital commutative algebra on the linear space $\Cal T$ spanned by nonempty non-planar rooted trees. We list all rooted trees up to degree $5$: 
$$
	\racine,\hskip 8mm  \arbrea,\hskip 8mm  \arbreba,\  \arbrebb,\hskip 8mm  \arbreca,\ \arbrecb,\ \arbrecc,\ \arbrecd,\hskip 8mm 
	 \arbreda,\ \arbredb,\ \arbredc,\ \arbredd,\ \arbrede,\ \arbredf,\
     \arbredz,\ \arbredg,\ \arbredh \ldots .
$$ 
The empty set is denoted $\un$, and is the unit. A {\sl rooted forest\/} is a finite collection $s=(t_1,\ldots, t_n)$ of rooted trees, which we simply denote by the (commutative) product $t_1 \cdots t_n$. Recall that the operator $B_+$ associates to the forest $s$ the tree $B_+(s)$ obtained by grafting the connected components $t_j$ on a common new root. $B_+(\un)$ is the unique rooted tree $\begin{matrix}  \\[-0.5cm]  \racine \end{matrix}$ with only one vertex. The Butcher--Connes--Kreimer coproduct on a rooted tree $t$ is described in terms of admissible cuts as follows: 
 \begin{equation*}
    \Delta_{\makebox{{\tiny{BCK}}}}(t) = t \otimes \un + \un \otimes t + \sum_{c \in {\tiny{\mop{Adm}(t)}}} P^c (t)\otimes R^c(t).
 \end{equation*}
Here $\mop{Adm}(t)$ is understood as the set of {\it{admissible cuts}} of a tree, i.e. the set of collections of edges such that any path from the root to a leaf contains at most one edge of the collection\footnote{In order to make this picture completely correct, we must stress that for any nonempty tree two admissible cuts must be associated with the empty collection: the empty cut and the total cut.}. We denote by $P^c(t)$ (resp. $R^c(t)$) the pruning (resp. the trunk) of $t$, i.e. the subforest formed by the edges above the cut
$c \in \mop{Adm}(t)$ (resp. the subforest formed by the edges under the cut). Note that the trunk of a tree is a tree, but the pruning of a tree may be a forest. An {\sl elementary cut\/} is a cut of only one edge. See \cite{Hoffman, Holtkamp1, Holtkamp2, Manchon} for more details on the combinatorics of rooted trees and Hopf algebras.

%%%%%%%%%%%%%%%%%%%%%%%%%%%%%%%%%%%%%%%%%%%%%

\subsection{Two isomorphic Hopf algebras of rooted trees}
\label{ssect:2hopf}

%%%%%%%%%%%%%%%%%%%%%%%%%%%%%%%%%%%%%%%%%%%%%

\subsubsection{Groups associated with augmented operads}
\label{sssect:GrOp}

Following \cite{Chap1}, we introduce an {\sl augmented operad\/}, which is an operad $\Cal P$ such that $\mop{dim}\Cal P_0$ and $\mop{dim}\Cal P_1=1$, i.e.~such that there is no $0$-ary operation, and such that the only $1$-ary operation is the unit $e$. The group $G_{\Cal P}$ is defined in \cite{Chap1} as the group of invertible elements in the product:
$$
	\prod_{n\ge 1}(\Cal P_n)_{S_n},
$$
which is the completed free $\Cal P$-algebra with one generator. An element $g=(g_n)_{n\ge 1}$ in this product is invertible if and only if its first component $g_1$ is nonzero. We will consider a slightly smaller group:
\begin{equation*}
	G_{\Cal P}^e:=\{g=(g_n)_{n\ge 1},\,g_1=e\}.
\end{equation*}
The advantage of this definition is the pro-nilpotency property. The associated Lie algebra is given by:
\begin{equation*}
	\g g_{\Cal P}^e:=\{x=(x_n)_{n\ge 1},\,x_1=0\}, 
\end{equation*}
with the Lie bracket given by:
\begin{equation*}
	[x,y]:=(x\curvearrowleft y) - (y\curvearrowleft x)
		 =\sum_{m,n\ge 2}(x_m\curvearrowleft y_n)-(y_n\curvearrowleft x_m)
\end{equation*}
where, from an operadic point of view:
\begin{equation*}
	x_m\curvearrowleft y_n:={\sum_{i=1}^m}x_m\circ_i y_n
			={\sum_{i=1}^m}\gamma(x_m;\underbrace{e,\ldots ,e}_{i-1\smop{ terms}},y_n,e,\ldots ,e).
\end{equation*}
The operation $\curvearrowleft$ defined above is right pre-Lie \cite{ChaLiv}, i.e.~we have:
\begin{equation}
\label{rightpL}
	(x\curvearrowleft y)\curvearrowleft z-x\curvearrowleft (y\curvearrowleft z)
		=(x\curvearrowleft z)\curvearrowleft y-x\curvearrowleft (z\curvearrowleft y).
\end{equation}
Of course, we also could consider the left pre-Lie operation $\curvearrowright$ defined by $x \curvearrowright y:=y\curvearrowleft x$, subject to the left pre-Lie relation:
\begin{equation}
\label{pre-Lie}
	(x \curvearrowright y) \curvearrowright z-x \curvearrowright (y \curvearrowright z) 
		= (y \curvearrowright x) \curvearrowright z-y \curvearrowright (x \curvearrowright z).
\end{equation}
The reader immediately verifies that $[x,y]=x\curvearrowleft y-y\curvearrowleft x=-(x \curvearrowright y-y \curvearrowright x)$.

%%%%%%%%%%%%%%%%%%%%%%%%%%%%%%%%%%%%%%%%%%%%%

\subsection{A Hopf algebra structure on reduced planar forests}
\label{ssect:Hopfpbt}

We now define a graded connected Hopf algebra structure on planar reduced forests, with grading given by the total number of inner vertices. First, we extend $\Cal T^{red}_{pl}$ to the free noncommutative algebra of reduced planar rooted forests, denoted by $\Cal H^{red}_{pl}$, with the one-edge tree $|$~as unit and the multiplication given by concatenation. We define a coproduct on reduced planar trees in terms of {\sl admissible cuts\/} of a tree $t \in \Cal T^{red}_{pl}$, i.e.~a (possibly empty) subset $c$ of edges {\sl not connected to a leaf\/} with the rule that along any path from the root of $t$ to any of its leaves there is at most one edge in $c$. The edges in $c$ are naturally ordered from left to right. To any admissible cut $c$ always corresponds then a unique subforest $P^c(t)$, the {\sl pruning\/}, obtained by concatenation of the subtrees obtained by cutting the edges in $c$, in the order defined as above. Then we define the coproduct:
$$
	\Delta_2(t)=\sum_{c\in\smop{Adm }t}P^c(t)\otimes R^c(t),
$$    
where $R^c(t)$ is the trunk, obtained by replacing each subtree of $P^c(t)$ with a single leaf. Note that the trunk of a tree is a tree, but the pruning of a tree may be a forest. We present a few examples:
 \allowdisplaybreaks{
\begin{eqnarray*}
	\Delta_2(\tree\hspace{0.1cm} ) &=& \tree\ \otimes |\ + | \otimes\! \tree \\[0.2cm]
	\Delta_2(\scalebox{0.8}{\treeA}\hspace{0.1cm} ) &=& \scalebox{0.8}{\treeA}\otimes |\ + | \otimes\! \scalebox{0.8}{\treeA}
										\hspace{0.3cm} + \  \tree \otimes\! \tree \\[0.2cm]
	\Delta_2(\scalebox{0.8}{\treeB}\hspace{0.1cm} ) &=& \scalebox{0.8}{\treeB} \otimes |\ + | \otimes\! \scalebox{0.8}{\treeB}
										\hspace{0.3cm} + \  \tree \otimes\! \tree \\[0.2cm]
	\Delta_2(\scalebox{0.8}{\treeC}\hspace{0.1cm} ) &=& \scalebox{0.8}{\treeC} \otimes |\ + | \otimes\! \scalebox{0.8}{\treeC}
										\hspace{0.3cm} + \  \scalebox{0.8}{\treeA} \otimes\! \tree 
		   									     \hspace{0.3cm} + \  \tree \otimes\! \scalebox{0.8}{\treeA} \\[0.2cm]
	\Delta_2(\scalebox{0.8}{\treeD}\hspace{0.1cm} ) &=& \scalebox{0.8}{\treeD} \otimes |\ + | \otimes\! \scalebox{0.8}{\treeD}
										 	      \hspace{0.3cm} + \  \tree \otimes\! \scalebox{0.8}{\treeA}
											         \hspace{0.3cm} + \  \tree \otimes\! \scalebox{0.8}{\treeB}
										 \hspace{0.3cm} + \  \tree\; \tree \otimes\! \tree \\[0.2cm]
	\Delta_2\scalebox{0.8}{(\treeE}\hspace{0.1cm} ) &=&\scalebox{0.8}{\treeE} \otimes |\ + | \otimes\! \scalebox{0.8}{\treeE}
										\hspace{0.3cm} + \  \scalebox{0.8}{\treeB} \otimes\! \tree 
		   									     \hspace{0.3cm} + \  \tree \otimes\! \scalebox{0.8}{\treeB} \\[0.2cm]
	\Delta_2(\scalebox{0.8}{\treeF}\hspace{0.1cm} ) &=& \scalebox{0.8}{\treeF}  \otimes |\ + | \otimes\! \scalebox{0.8}{\treeF}
										\hspace{0.3cm} + \  \scalebox{0.8}{\treeA} \otimes\! \tree 
		   									     \hspace{0.3cm} + \  \tree \otimes\! \scalebox{0.8}{\treeB} \\[0.2cm]
	\Delta_2(\scalebox{0.8}{\treeG}\hspace{0.1cm} ) &=& \scalebox{0.8}{\treeG}\otimes |\ + | \otimes\! \scalebox{0.8}{\treeG}
										\hspace{0.3cm} + \  \scalebox{0.8}{\treeB} \otimes\! \tree 
		   									     \hspace{0.3cm} + \  \tree \otimes\! \scalebox{0.8}{\treeA} \\[-.2cm]
	\Delta_2(\scalebox{0.8}{\treebc}) &=& \scalebox{0.8}{\treebc} \otimes | + | \otimes\! \scalebox{0.8}{\treebc}
										 + \  \scalebox{0.8}{\treec}   \hspace{0.3cm}\ \otimes\! \tree					
\end{eqnarray*}}
										 
We remark here that several Hopf subalgebras are readily identified. First, the binary forests obviously form a Hopf subalgebra $\Cal H_{pl}^{bin}$ of $\Cal H_{pl}^{red}$, which in turn contains two other Hopf subalgebras, i.e.~the Hopf subalgebra $\Cal H_{r, pl}^{bin}\subset \Cal H_{pl}^{bin}$ (resp. $\Cal H_{l, pl}^{bin} \subset \Cal H_{pl}^{bin}$) of {\sl{right- (resp. left)-combed}} binary planar rooted forests, generated by the trees $t^{(n)}_r$ (resp. $t^{(n)}_l$) recursively defined by $t^{(1)}_{r}:= | =:t^{(1)}_{l}$ and $t^{(n)}_r:= | \vee t^{(n-1)}_r$, (resp. $t^{(n)}_l:= t^{(n-1)}_l \vee |$), $n>1$. Also, observe that the trees with only one inner vertex, let us call them {\sl{reduced corollas}},  are all primitive.

It is immediate to adapt this construction to the with setting decoration described in Paragraph \ref{sect:deco}. Details are left to the reader.

%%%%%%%%%%%%%%%%%%%%%%%%%%%%%%%%%%%%%%%%%%%%%

\subsection{The associated pre-Lie structure}
\label{ssect:prelie}

Let $(\Cal H_{pl}^{red})^\circ$ be the graded dual of $\Cal H_{pl}^{red}$. We consider the normalized dual basis $(\delta'_t)$ of the basis of forests, defined by:
\begin{equation*}
	<\delta'_{t_1\cdots t_k},\,t_1\cdots t_k>=\sigma(t_1)\cdots\sigma(t_k)
\end{equation*}
where $\sigma(t_j)$ is the symmetry factor of the tree $t_j$, and $<\delta'_{t_1\cdots t_k},\,s>=0$ if $s$ is a forest different from $t_1\cdots t_k$. The correspondence $t_1\cdots t_k\mapsto\delta'_{t_1\cdots t_k}$ yields a linear isomorphism $\delta':\Cal H_{pl}^{red}\to (\Cal H_{pl}^{red})^\circ$. If $t$ and $u$ are planar reduced trees, $\delta'_t$ and $\delta'_u$ are infinitesimal characters of the Hopf algebra $\Cal H_{pl}^{red}$, hence so is the Lie bracket defined in terms of the convolution product, $[\delta'_t,\,\delta'_u]=\delta'_t \star_2 \delta'_u-\delta'_u \star_2 \delta'_t$. Recall that an infinitesimal character maps the one-edge tree $|$ as well as any forest $t_1\cdots t_k$, $k>1$ to zero. The definition of the convolution product yields:
\begin{equation*}
	[\delta'_t,\,\delta'_u]=\delta'_{t{\curvearrowright_\sigma} u-u{\curvearrowright_\sigma} t},
\end{equation*}
where we define for any reduced planar rooted tree $t$:
\begin{equation*}
	t{\curvearrowright_\sigma} u=\sum_{v\in T_{pl}^{red}}\frac{\sigma(t)\sigma(u)}{\sigma(v)}N(t,u,v)v.
\end{equation*}
The coefficient $N(t,u,v)$ is the number of elementary cuts $c$ of the tree $v$ (in the sense of
the previous subsection) such that $P^c(v)=t$ and $R^c(v)=u$. The coefficient:
\begin{equation*}
	M(t,u,v):=\frac{\sigma(t)\sigma(u)}{\sigma(v)}N(t,u,v)
\end{equation*}
is the number of ways to graft $t$ on a leaf of the tree $u$ in order to obtain the tree $v$. Hence $t {\curvearrowright_\sigma} u$ is the sum of all the possible graftings of $t$ on $u$. It is well-known that the left pre-Lie relation holds:
\begin{equation*}
	s{\curvearrowright_\sigma}(t{\curvearrowright_\sigma} u)-(s{\curvearrowright_\sigma} t){\curvearrowright_\sigma} u
		=t{\curvearrowright_\sigma}(s{\curvearrowright_\sigma} u)-(t{\curvearrowright_\sigma} s){\curvearrowright_\sigma} u.
\end{equation*}
Namely both sides are expressed as the sum of all possible ways of grafting $s$ and $t$ on two different leaves of $u$. The associated Lie algebra structure on ${\Cal T}^{red}_{pl}$ is defined by $[t,u]:=t{\curvearrowright_\sigma} u-u{\curvearrowright_\sigma} t$, and gives rise to the Lie algebra of the (pro-nilpotent) group of characters, that is, multiplicative maps on the Hopf algebra $\Cal H^{red}_{pl}$, which identifies with the group $G_{\Cal F(V)}^{e,\smop{op}}$ associated with the free operad $\Cal F(V)$ on the $\mathbb{S}$-object $V=(k[S_n])_{n\ge 1}$, but with multiplication reversed\footnote{This is due to the fact that the Lie algebra structure comes from a left pre-Lie operation.}. Let us remark that the commutative Hopf algebra, which follows via the Cartier--Milnor--Moore theorem from the group $G_{\Cal F(V)}^{e,\smop{op}}$, is not isomorphic to $\Cal H^{red}_{pl}$, but is just a quotient.

\smallskip 

The same construction with planar binary trees yields the group of characters of the Hopf algebra $\Cal H^{bin}_{pl}$, which identifies with the group $G_{\Cal B}^{e,\smop{op}}$ associated with the free binary operad (with multiplication reversed). The free binary operad is the free operad $\Cal F(W)$ on the $\mathbb{S}$-object $W$ such that $W_1=k$, $W_2=k[S_2]$ and $W_n=\{0\}$ for $n\ge 3$. Finally the same construction with $\Cal I$-decorated trees (with the notations of Paragraph \ref{sect:deco}) yields the group $G^{e,\smop{op}}_{\Cal F(V_{\Cal I})}$, where $\Cal F(V_{\Cal I})$ is the free operad on the $\mathbb{S}$-object $V_{\Cal I}$ defined by $(V_{\Cal I})_1=k$ and $(V_{\Cal I})_n=k^{|I_n|}\otimes k[S_n]$ for $n\ge 2$.

%%%%%%%%%%%%%%%%%%%%%%%%%%%%%%%%%%%%%%%%%%%%%

\subsection{A Hopf algebra structure on planar rooted hyperforests}
\label{ssect:Hopfhypertrees}

We extend the linear isomorphism $\Phi:\Cal T_{pl}^{red}\to H\Cal T_{pl}$ to a graded algebra isomorphism still denoted by $\Phi:\Cal H_{pl}^{red}\to H\Cal H_{pl}$, where $H\Cal H_{pl}$ stands for the free noncommutative algebra of rooted planar hyperforests. The grading is given by the total number of {\sl edges\/}. The Hopf algebra structure on $\Cal H_{pl}^{red}$ can be transferred on $H\Cal H_{pl}$ by $\Phi$. The coproduct $\Delta=(\Phi\otimes\Phi)\circ\Delta_2\circ \Phi^{-1}$ can then be made explicit as  follows.

\medskip

We introduce the concept of {\sl right admissible cut\/} in the spirit of Munthe-Kaas and Wright \cite{MKW}. For any vertex $v \in V(t)$ we denote by $f(v)$ its fertility, i.e.~the number of edges with root $v$. Recall that we work with planar hypertrees. Hence, we may enumerate the incoming edges of each vertex $v \in V(t)$ counterclockwise from $1$ to $f(v)$. For any vertex $v$ and for any $i \in \{1,\ldots ,f(v)\}$ the  {\sl i$^{\smop{th}}$ single right vertex-cut\/} associated to $v$ is the subset $c^{(i)}_v \subset E(t)$ of the $i$ first edges with root $v$ with respect to the order above. To each single right vertex-cut $c^{(i)}_v$ we may associate a sub-hypertree $P^{c^{(i)}_v}(t)$ obtained from $t$ by removing the edges $c^{(i)}_v(t)$ from the vertex $v$ in $t$ and grafting them to a new root resulting in a single planar rooted hypertree. We denote by $R^{c^{(i)}_v}(t)$ the remaining tree. A {\sl right vertex-cut \/} $C$ is a (possibly empty) collection of single right vertex-cuts.  A (right) vertex-cut $C$ is called {\sl admissible\/} if any path from the root to any vertex of $t$ encounters at most one single right vertex-cut. The single vertex-cuts in an admissible $C$ are naturally ordered from left to right, thus giving rise to a planar hyperforest $P^C(t)$. We denote by $R^{C}(t)$ the remaining tree. By ${\rm{RAdm}}(t)$ we denote the set of admissible right vertex-cuts. We define in terms of  admissible (right) vertex-cuts the following coproduct:
\begin{equation*}
	\Delta(t)=\sum_{C \in {\rm{RAdm}}(t)} P^{C}(t) \otimes R^{C}(t),
\end{equation*}
\goodbreak
We list a few coproducts below. Observe the conservation of the number of edges. 
 \allowdisplaybreaks{	
\begin{eqnarray*}
\Delta(\!\racine)&=&\racine\otimes \racine, \quad\ 
\Delta(\arbrea)= \arbrea \otimes \racine + \racine \otimes  \arbrea\\
&&\\
\Delta(\arbreba)&=& \arbreba \otimes \racine + \racine \otimes  \arbreba  + \arbrea \otimes \arbrea, \quad\
\Delta(\arbrebb)= \arbrebb \otimes \racine + \racine \otimes  \arbrebb  + \arbrea \otimes \arbrea\\
&&\\
\Delta(\arbreca)&=&\arbreca \otimes \racine + \racine \otimes \arbreca
					+ \arbrea \otimes \arbreba + \arbreba \otimes \arbrea, \quad\
\Delta(\arbrecb)=\arbrecb \otimes \racine + \racine \otimes  \arbrecb
					+\arbrea \otimes \arbreba
					+\arbrebb \otimes \arbrea\\
\Delta(\arbrecc)&=&\arbrecc \otimes \racine + \racine \otimes  \arbrecc
					+\arbrea \otimes \arbreba
											+\arbrea \arbrea\otimes \arbrea
					+\arbrea \otimes \arbrebb\\
\Delta(\arbrece)&=&\arbrece \otimes \racine + \racine \otimes  \arbrece
					+\arbrea \otimes \arbrebb
					+\arbreba \otimes \arbrea, \quad\
\Delta(\arbrecd)=\arbrecd \otimes \racine + \racine \otimes  \arbrecd
					+\arbrea \otimes \arbrebb
					+\arbrebb \otimes \arbrea\\
&&\\
\Delta(\arbreda)&=&\arbreda \otimes \racine + \racine \otimes \arbreda  
					+ \arbrea \otimes \arbreca
					+ \arbreba \otimes \arbreba
					+ \arbreca \otimes \arbrea\\
\Delta(\arbredb)&=&\arbredb \otimes \racine + \racine \otimes \arbredb  
					+ \arbrea \otimes \arbreca
					      + \arbrebb \otimes \arbreba
					+ \arbrecb \otimes \arbrea\\
\Delta(\arbredf)&=&\arbredf \otimes \racine + \racine \otimes \arbredf 
					+ \arbrea \otimes \arbrece
				    	      + \arbrea \otimes \arbrecc
					+\arbrea \arbrea\otimes \arbrebb
						+\arbreba \otimes \arbreba\\
\Delta(\arbredg)&=&\arbredg \otimes \racine + \racine \otimes \arbredg
					      + \arbrea \otimes \arbrecd
					      + \arbrea\arbrea \otimes \arbrebb
					      + \arbrea \otimes \arbrece
				    	      + \arbrecc\otimes \arbrea \\		
\Delta(\arbrede)&=&\arbrede \otimes \racine + \racine \otimes \arbrede
					      + \arbrea \otimes \arbreca
					      + \arbrea\arbrea \otimes \arbreba
					      	     + \arbrea\arbreba \otimes \arbrea
					      + \arbrea \otimes \arbrecc
				    	      + \arbreba\otimes \arbrebb\\
&&\\
\Delta(\scalebox{0.55}{\haced})&=&\scalebox{0.55}{\haced}\otimes\racine
									+\racine\otimes\scalebox{0.55}{\haced} 
									+\scalebox{0.55}{\hab}\otimes\arbrea
\end{eqnarray*}}

We note that via the bijection $\Phi$ we identify the Hopf subalgebras $\Phi(\Cal H_{pl}^{bin} )=\Cal H_{pl} \subset H\Cal H_{pl}$ and  $\Phi(\Cal H_{l,pl}^{bin} )=\Cal H_{pl}^{lad} \subset \Cal H_{pl}$  and  $\Phi(\Cal H_{r,pl}^{bin} )=\Cal H_{pl}^{cor} \subset \Cal H_{pl}$ of ladder trees and corollas, respectively. Reduced corollas with $n$ leaves are mapped to blobs with $n$ vertices drawn on the circle delimiting the blob.

This Hopf algebra structure is related to the pre-Lie structure $\surj$ defined on $H\Cal T_{pl}$ by:
\begin{equation*}
	s_1\surj s_2 := \Phi\big(\Phi^{-1}(s_1){\curvearrowright_\sigma}\Phi^{-1}(s_2)\big),
\end{equation*}
where ${\curvearrowright_\sigma}$ is the pre-Lie product defined in Subsection \ref{ssect:prelie}. The associated Lie algebra is of course isomorphic to the one defined in same subsection. It is then another presentation of the opposite Lie algebra of the pro-nilpotent group $G_{\Cal F(V)}^e$ associated with the free operad on the $\mathbb{S}$-object $V$ defined in Paragraph \ref{ssect:prelie}. The same construction with planar rooted trees gives back the group $G_{\Cal B}^e$ associated with the free binary operad (modulo reversing the multiplication or, what is the same, changing the sign of the Lie bracket). The same construction holds for $\Cal I$-decorated hypertrees, leading to another presentation of the opposite Lie algebra of the pro-nilpotent group $G_{\Cal F_{V(\Cal I)}}^e$ associated with the free operad on the $\mathbb{S}$-object $V(\Cal I)$: details are left to the reader.\\

\bigskip

{\bf{Acknowledgements}} We would like to thank F.~Chapoton for very helpful discussions. We thank H.~Munthe-Kaas and A.~Lundervold for comments. The first author is supported by a Ram\'on y Cajal research grant from the Spanish government. We thank the GDR Renormalisation for support.\\

%%%%%%%%%%%%%%%%%%%%%%%%%%%%%%%%%%%%%%%%%%%%%
%%%%%%%%%%%%%%%%%%%%%%%%%%%%%%%%%%%%%%%%%%%%%
%%%%%%%%%%%%%%%%%%%%%%%%%%%%%%%%%%%%%%%%%%%%%

\end{document}